\documentclass{amsart}
\usepackage{amsthm}
\usepackage{graphicx}
\usepackage{amssymb}
\usepackage{psfrag}
\title{$6j$-symbols, hyperbolic structures and the Volume Conjecture}

\author[F. Costantino]{Francesco Costantino}
\address{Institut de Recherche Math\'ematique Avanc\'ee\\
  Rue Ren\'e Descartes 7\\
  67084 Strasbourg, France}
\email{costanti@math.u-strasbg.fr}
\newtheorem{teo}{Theorem}[section]
\newtheorem{lemma}[teo]{Lemma}
\newtheorem{prop}[teo]{Proposition}

\newtheorem{conj}[teo]{Conjecture}
\newtheorem{question}[teo]{Question}

\theoremstyle{definition}
\newtheorem{defi}[teo]{Definition}
\newtheorem{rem}[teo]{Remark}
\newtheorem{example}[teo]{Example}

\theoremstyle{remark}
\newtheorem{prof}[teo]{Proof of}

\def\mn{\mathbb{N}}

\def\mc{\mathbb{C}}
\def\mr{\mathbb{R}}
\def\mz{\mathbb{Z}}
\def\mq{\mathbb{Q}}
\def\mh{\mathbb{H}}

\def\ns{\negmedspace}
\def\nns{\negthickspace}

\begin{document}
%\dedicatory{This is a pre-preprint.  Please provide us with comments!}

\begin{abstract}
We compute the asymptotical growth rate of a large family of $U_q(sl_2)$ $6j$-symbols and we interpret our results in geometric terms by relating them to volumes of hyperbolic truncated tetrahedra. We address a question which is strictly related with S.Gukov's generalized volume conjecture and deals with the case of hyperbolic links in connected sums of $S^2\times S^1$. We answer this question for the infinite family of fundamental shadow links.
\end{abstract}
%\begin{asciiabstract}
%We compute the asymptotical growth rate of a large family of $U_q(sl_2)$ $6j$-symbols and we interpret our results in geometric terms by relating them to volumes of hyperbolic truncated tetrahedra. We propose an extension of S. Gukov's generalized volume conjecture to cover the case of hyperbolic links in $S^3$ or connected sums of $S^2\times S^1$. We prove this conjecture for the infinite family of fundamental shadow links.   
%\end{asciiabstract}

\maketitle

\tableofcontents
\section{Introduction}
The study of $6j$-symbols of $U_q(sl_2)$ plays a central role in mathematical physics (for the study of scattering of particles), in algebra (for the study of representation theory of $U_q(sl_2)$), in 3-dimensional geometry (for the computation of volumes of tetrahedra in curved spaces), and in quantum topology (for the computation of quantum invariants). For our purposes, we will define a $6j$-symbol as an explicit rational function of the variable $q$ which is parametrized by $6$ half-integer numbers.
The present paper is devoted to give a partial answer to the following question (see Definition \ref{def:admissibility} for the definition of  ``admissible" $6$-tuple):
 \begin{question}[Asymptotical behavior of $6j$-symbols]\label{que:main}
Given $\vec \theta\in \mr_+^6$ and a sequence $\vec b^n=(b^n_0,b^n_1,b^n_2,b^n_3,b^n_4,b^n_5)\in \mn^6,\ n\in \mn$ (or, more in general of admissible $6$-tuples of half integers) such that $\lim_{n\to \infty} \frac{b^n_i}{n}=\theta_i$.
What is the asymptotical behavior of the lowest order term of the Laurent series expansion around $q=Exp(\frac{2\pi \sqrt{-1}}{n})$ of $$ \setlength\arraycolsep{1pt}
\left( \begin{array}{ccc}
b^n_0 & b^n_1 & b^n_2\\
b^n_3  & b^n_4 & b^n_5\\
\end{array}\right)$$ when $n\to \infty$?
\end{question}

An answer to Question \ref{que:main} has been provided by Y.Taylor and C.Woodward (\cite{TW}) in the case when $\vec \theta$ satisfies a set compatibility conditions we called ``Reshetikhin-Turaev'' conditions. Their result shows that, in this case, the asymptotical behavior is governed by the geometry of a spherical tetrahedron whose edge lengths are proportional to the entries of $\vec \theta$ and so, in particular, not by an hyperbolic object.

In the present paper we study a completely disjoint set of cases: we define different compatibility conditions (see Definition \ref{def:radmissibility}) identifying a set of $\vec\theta's$ which we called \emph{hyperbolic} and study Question \ref{que:main} for this set. 
The $\theta's$ we consider are such that there exists a unique hyperbolic hyperideal tetrahedron whose internal dihedral angles are $\alpha_i=|2\pi(\theta_i-\frac{1}{2})|$. Recall that Murakamy and Yano (\cite{MYa}) provided a symmetric, real analytic formula for the volume of a proper (i.e.  entirely contained in $\mathbb{H}^3$) hyperbolic tetrahedron; later  A. Ushijima (\cite{Us}) proved that if one applies the same formula to the set of angles of an hyperideal hyperbolic tetrahedron, one gets the volume of the truncation of tetrahedron by means of geodesic planes orthogonal to its faces. 
The following (Theorem \ref{teo:asymptotics} below) gives an answer to a question asked by Taylor and Woodward (Question c, \cite{TW}):
\begin{teo}\label{teo:asymptoticsem}
Let $\vec \theta\in \mr_+^6$ be a hyperbolic $6$-tuple, and let $\vec b^n=(b_0^n,b_1^n,b_2^n,b_3^n,b_4^n,b_5^n)\in \frac{\mn}{2}^6, n\in \mn$ be a sequence of admissible $6$-tuples of half integers such that $\lim_{n\to \infty} \frac{b^n_i}{n}=\theta_i$. It holds:
$$\lim_{n\to \infty} \frac{2\pi}{n}Log(|ev_n(\left( \begin{array}{ccc}
b^n_0 & b^n_1 & b^n_2\\
b^n_3  & b^n_4 & b^n_5\\
\end{array}\right)|)=2Vol(T(\vec\alpha))$$
where $Vol$ is the hyperbolic volume and $T(\alpha)$ is the hyperbolic truncated tetrahedron whose dihedral angles are $\alpha_i=|2\pi(\theta_i-\frac{1}{2})|$ and by $ev_n$ we denote the evaluation as explained later in the notation.
\end{teo}

The main application of the above result is related to a question strictly related with the Volume Conjecture.
\begin{conj}[Kashaev~\cite{K}, Murakamy-Murakamy~\cite{MM}]\label{conj:kmm}
Let $L$ be a link in $S^3$, and $J'_n(q)$ its $n^{th}$-Jones polynomial normalized so that $J'_n(\rm{unknot})=1, \forall n\geq 2$.
The following holds:
$$\lim_{n\to \infty} \frac{2\pi}{n} Log(|J'_n(e^{\frac{2\pi i}{n}})|)=Vol(S^3\setminus L)$$
\end{conj}
Before stating our question, let us recall some facts about the conjecture.
The VC in the above form has been formally checked for the Figure Eight knot (\cite{MM}), the Borromean link, for torus knots (\cite{KT}, \cite{MM3}), their Whitehead doubles (\cite{Zh}) and for the family of ``Whitehead chains" (\cite{Va}). Moreover, there is experimental evidence of its validity for the knots $6_3$, $8_9$ and $8_{20}$ (\cite{MMOTY}).

More recently, Baseilhac and Benedetti (\cite{BB1}, \cite{BB2}, \cite{BB3}) formulated different instances of general Volume Conjectures in the framework of their ``quantum hyperbolic field theory" (QHFT); however, Kashaev's VC is {\it not} a specialization of such QHFT Volume Conjectures.

In \cite{Cocv}, we proposed an ``extension" of the VC to include the case of links in $S^3$ or in connected sums of copies of $S^2\times S^1$: we defined colored Jones invariants taking values in the ring of rational functions for links in these more general ambient manifolds and proved the conjecture for the infinite family of ``fundamental shadow links"  (called universal hyperbolic links in the paper). These links were already studied in \cite{CT2} because of their remarkable geometric properties.
We stress that, in the case of links in $S^3$, our colored Jones invariants reduce to colored Jones polynomials and our extension of the Volume Conjecture reduces to Kashaev-Murakamy-Murakamy's conjecture.

In \cite{Gu}, S. Gukov proposed a generalization of the VC for the case of knots in $S^3$ to include evaluations of $J'_n$ in points near $e^{\frac{2\pi i}{n}}$ (see Section \ref{sec:gvc} for a precise statement). Gukov's Conjecture was proved only for the Figure Eight knot and for all torus knots (see \cite{MY}).
In the present paper, we address the following question which is related to a possible extension of Gukov's conjecture for hyperbolic links in $N=\# k <\nns S^2\times S^1\nns>$.
\begin{question}\label{conj:genforlinkssemplice}
Is it true that there exists a neighborhood $U$ of $\vec 0\in \mr^r$ such thatÊ for each $\vec a\in U\cap \mr^r$ and for each sequence $\vec b^n\in \mn^r$, such that $\lim_{n\to \infty}\frac{\vec b^n}{n}=\frac{\vec 1+\vec a}{2}$, the following holds?
$$\lim_{n\to \infty} \frac{2\pi}{n}Log(|ev_{n}(J_{\vec b^n})|) =Vol((N\setminus L)_a)$$
In the above expression, $J_{\vec b^n}$ is the colored Jones invariants of $L$ associated to the coloring $\vec b^n$ of $L$, $\vec 1\in \mr^r$ is the vector whose entries are all $1$, $(N\setminus L)_a$ is the hyperbolic structure whose holonomy around the meridian of $L_i$ is a matrix conjugated to an upper triangular matrix whose eigenvalues are $e^{\pm \pi i a_i }$. Ê
\end{question}
I am indebted with by Stavros Garoufalidis who pointed out to me that the answer to the above question is ``no" if $L$ is a knot in $S^3$ (see Section \ref{sec:gvc} for more details on this). On contrast, when the ambient manifold $N$ is not $S^3$, we shall exhibit infinitely many examples where the answer is positive. To do this, in Section \ref{sec:tetrahedra} we will recall the definition of the fundamental shadow links and prove a decomposition result for the hyperbolic structures on their complements (Proposition \ref{prop:kerchoff}).  
Then we will prove the following (Theorem \ref{teo:gcv}):
\begin{teo}
The answer to Question \ref{conj:genforlinkssemplice} is ``yes" for the infinite family of fundamental shadow links.
\end{teo}

{\bf Acknowledgments.} I wish to warmly thank Stephane Baseilhac, Riccardo Benedetti, Dylan Thurston and Vladimir Turaev for the helpful critics and suggestions they provided to me. During the development of the present research, I was supported by a Marie-Curie Intraeuropean Fellowship at IRMA (Strasbourg). I also benefited of a research period at Institut Fourier de Grenoble funded by ANR. I am indebted with Stavros Garoufalidis who, through his comments and suggestions, helped me to correct and improve the previous version of this paper.

{\bf Notation.} 
All the dihedral angles will be internal unless explicitly stated the contrary.
If $f:\mc\to \mc\cup\{ \infty\}$ is a meromorphic function, the \emph{evaluation of $f$ at $x\in \mc$}, denoted $ev_x(f)$, is the leading coefficient of the Laurent series expansion of $f$ at $x$. We will denote $ev_n(f)$ the evaluation of $f$ at $q_n=e^{\frac{2\pi i}{n}}$.
Note that if $f$ and $g$ are meromorphic functions then $ev_z(fg)=ev_z(f)ev_z(g)$ regardless the presence of poles or zeros at $z$.
\section{Truncated tetrahedra and fundamental shadow links}\label{sec:tetrahedra}
In this section we recall the definition of truncated tetrahedra and a result classifying them up to isometry. Then we define fundamental shadow links and show how to decompose their complements by means of truncated tetrahedra.
\subsection{Truncated tetrahedra and $D$-blocks}
Let $\Pi_1,\ldots \Pi_4$ be half spaces in $\mh^3$ having non-empty intersection and, through the embedding of the Klein model in $\mathbb{RP}^3$, let $H_i$ be the hyperplane delimiting $\Pi_i$. Let us suppose that $\mh^3\cap H_i\cap H_j\neq\emptyset,\ \forall i\neq j$ and that $v_t\doteq H_i\cap H_j\cap H_k\in \mathbb{RP}^3\setminus \mh^3$  for all distinct $i,j,k,t\in \{1,2,3,4\}$. Let also $K_i$ be hyperplanes in $\mh^3$ such that $K_i$ intersects $H_j,H_t,H_k$ orthogonally for all distinct $i,j,k,t\in \{1,2,3,4\}$ (one can check that $K_i$ indeed exists using the hyperboloid model for $\mh^3$) and let us suppose that $\mh^3\cap K_i\cap K_j=\emptyset, \forall i\neq j$.
\begin{defi}
An hyperbolic \emph{truncated tetrahedron} is the compact polyhedron in $\mh^3$ bounded by $H_1,\cup \ldots\cup H_4\cup K_1,\cup\ldots\cup K_4$. The \emph{faces} of the truncated tetrahedron $T$ are the hyperbolic polygons $T\cap H_i,\ i=1,\ldots 4$, its \emph{truncation faces} are the hyperbolic triangles $T\cap K_i,\ i=1\ldots 4$, its \emph{egdes} $e_{ij}$ are the geodesic arcs $T\cap H_i\cap H_j,\ \forall i\neq j$ and its \emph{angles} $\alpha_{ij}$ are the internal dihedral angles between $H_i$ and $H_j$ for all distinct $i$ and $j$. As limit cases, we allow faces $H_i$ and $H_j$ to be tangent in a point in $\partial \mh^3$: in this case $T$ is not compact (it has some ideal vertex), the edge $e_{ij}$ is just the point at infinity $H_i\cap H_j$ and has length $0$ and the dihedral angle $\alpha_{ij}$ is $0$.
\end{defi}
\begin{example}
The most remarkable example is when all the dihedral angles are $0$: in that case $T$ is a regular ideal octahedron and all the faces and truncation faces are ideal triangles intersecting orthogonally.
\end{example}

In the following we will use the following parametrization by dihedral angles of truncated tetrahedra (see for instance Theorem 2.2 of \cite{FP}):
\begin{teo}\label{teo:frigpetr}
A truncated tetrahedron with dihedral angles given by the entries of $\vec \alpha\in [0,\pi)^6$ exists and is unique up to isometry in $\mh^3$ iff for each truncation face the sum of the dihedral angles on the three edges touching the face is less than $\pi$.
\end{teo}

\begin{figure}
\psfrag{G}{\large $G$}
 \centerline{\includegraphics[width=6.5cm]{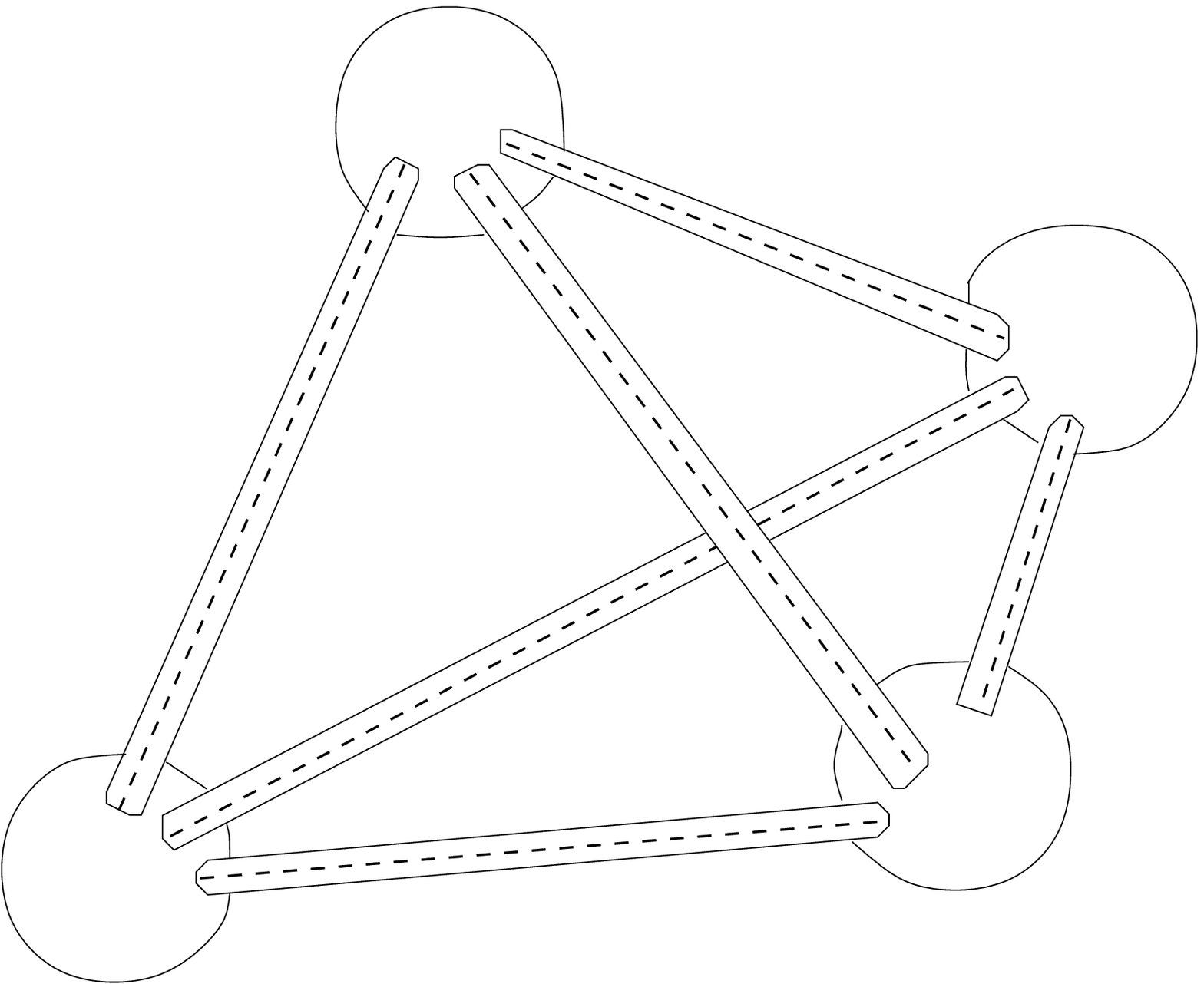}}
  \caption{}
  \label{fig:tetrahedron}
\end{figure}

\begin{defi}[D-blocks]
Let $T(\vec\alpha)$ be a truncated tetrahedron with dihedral angles $\vec \alpha\in [0,\pi)^6$. If $\vec u\in [0,2\pi)^6$ we call the $D$-block $D(\vec u)$ the hyperbolic $3$-manifold obtained by gluing two copies $T^+(\frac{\vec u}{2})$ and $T^-(\frac{\vec u}{2})$ of $T(\frac{\vec u}{2})$ oriented differently through the identity map on their faces. 
\end{defi}
When $\vec u\in (0,2\pi)^6$, $D(\vec u)$ is topologically $S^3\setminus 4B^3$ with an hyperbolic structure having cone angle singularities around the arcs dotted in Figure \ref{fig:tetrahedron}. The total angle around the arc corresponding to the edge $e_{ij}$ of $T(\frac{\vec u}{2})$ is the corresponding component of $\vec u$.
When some component of $\vec u$ is $0$, the hyperbolic structure on $S^3\setminus 4B^3$ has an annular cusp around the corresponding arc.

In any case $\partial D(\vec u)$ is the union of $4$ totally geodesic spheres $S_i,\ i=1,\ldots 4$ (round in Figure \ref{fig:tetrahedron}) each formed by the two truncation faces $K^+_i$ and $K^-_i$ respectively of $T^+(\vec\alpha)$ and $T^-(\vec\alpha)$, where $\vec \alpha=\frac{\vec u}{2}$. Each $S_i$ has $3$ cone-angle singularities of total angle $2\alpha_{jk},2\alpha_{jt},2\alpha_{tk}$ where $i,j,k$ and $t$ are distinct and a cone-angle $0$ singularity is a cusp.
\subsection{Murakamy-Yano's formula}

Let $T(\vec\alpha)$ be an hyperbolic (proper or truncated) tetrahedron whose dihedral angles are the entries of $\vec\alpha\in [0,\pi)^6$ with $\alpha_i$ and $\alpha_{i+3}$ corresponding to opposite angles for $i\in\{0,1,2\}$. Let $A_i=Exp(\sqrt{-1}\alpha_i)$ and let 
\begin{gather*}
U(z,T)=\frac{1}{2}(Li_2(z)\ns+\ns Li_2(zA_0A_1A_3A_4)+\ns Li_2(zA_0A_2A_3A_5)\ns+\ns Li_2(zA_2A_1A_5A_4)\ns-\\ 
-\ns Li_2(-zA_0A_1A_2)\ns -\ns Li_2(-zA_0A_4A_5)\ns -\ns Li_2(-zA_3A_1A_5)\ns -\ns Li_2(-zA_0A_4A_2))
\end{gather*} 

where $Li_2(z)=-\int_0^z \frac{Log(1-t)}{t}dt$. Let $z_{\pm}$ the two non-trivial solutions to the equation: 
\begin{equation}\label{eqnz}
\frac{d}{dz}U(z,T)=\frac{\pi \sqrt{-1}}{z}k\ \ \ \  k\in \mz
\end{equation}
It turns out that whenever $T(\vec\alpha)$ is hyperbolic (proper or truncated) then the norms of $z_\pm$ are $1$ (see the proof of the Lemma in Section 1.1 of \cite{MY}). Let us call in what follows $z_+$ ($z_-$) the solution with positive (negative) imaginary part.
Let moreover be:

\begin{gather*}
\hat{\Delta}(a,b,c)=-\frac{1}{4}(Li_2(-abc^{-1})+Li_2(-bca^{-1})+Li_2(-acb^{-1})+Li_2(-a^{-1}b^{-1}c^{-1})\\ 
+(Log\ a)^2+(Log\ b)^2+(Log\ c)^2)\\
\hat{\Delta}(T)=\hat{\Delta}(A_0,A_1,A_2)+\hat{\Delta}(A_0,A_4,A_5)+\hat{\Delta}(A_3,A_1,A_5)+\hat{\Delta}(A_3,A_4,A_2)\\ 
+\frac{1}{2}(Log\ A_0\ Log\ A_3+ Log\ A_1\ Log\ A_4+Log\ A_2\ Log\ A_5)
\end{gather*}
\begin{equation}\label{murakami-yano}
V(T)=\frac{\sqrt{-1}}{2}(U(z_-,T)-U(z_+,T)+Log(z_+)z_+\frac{\partial U(z,T)}{\partial z}|_{z=z_+}-Log(z_-)z_-\frac{\partial U(z,T)}{\partial z}|_{z=z_-})
\end{equation}

\begin{teo}[Murakamy-Yano \cite{MY}, Ushijima \cite{Us}]\label{teo:MY}
\begin{equation}\label{MYformula}
Vol(T(\vec\alpha))=Im(U(z_+,T)+\hat{\Delta}(T))=-Im(U(z_-,T)+\hat{\Delta}(T))=-V(T)
\end{equation}
\end{teo}
A very nice geometric interpretation of the above formula has been provided by P. Doyle and G. Leibon (\cite{DL}, see also Y. Mohanty's paper \cite{Mo}).  
\subsection{Fundamental shadow links}

\begin{figure}
\psfrag{G}{\large $G$}
\psfrag{L}{\large $L$}
\psfrag{0}{\large $0$}
 \centerline{\includegraphics[width=14.5cm]{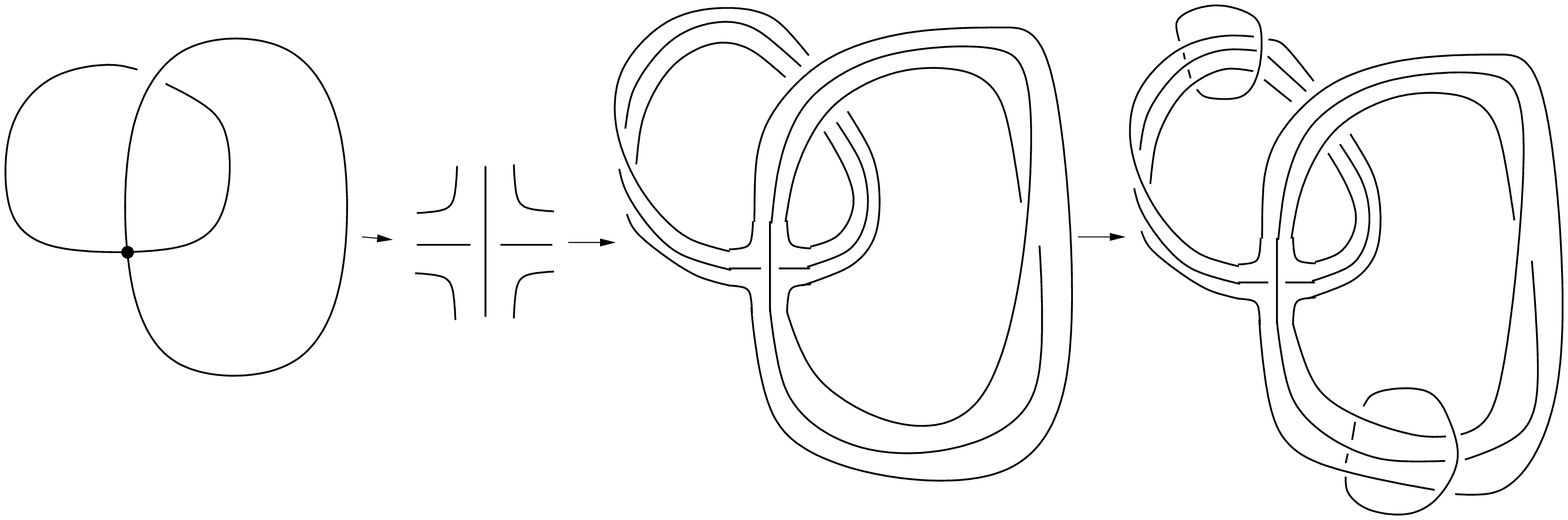}}
  \caption{}
  \label{fig:uhl}
\end{figure}
\begin{defi}\label{def:fundshlink}
A \emph{fundamental shadow link} is a link $L$ admitting a surgery presentation obtained as follows (see Figure \ref{fig:uhl}):
\begin{enumerate}
\item Consider a planar diagram of a $4$-valent graph $G$ embedded in $\mr^3$, and $T$ be a maximal subtree of $G$;
\item Replace each vertex of $G$ by a diagram composed of six strands as shown in Figure \ref{fig:uhl} and each edge by a $3$-braid;
\item Encircle each braid corresponding to an edge of $G\setminus T$ by $0$-framed meridians (in the figure $T=\emptyset$). 
\end{enumerate}
The number of vertices of $G$ is called the \emph{shadow-complexity} of $L$ and is denoted $sc(L)$.
\end{defi}

The following summarizes the properties of fundamental shadow links:
\begin{teo}[F.C. \& D.P. Thurston,\ \cite{CT2}]\label{teo:ioedylan}
Let $L$ be a fundamental shadow link obtained from a graph $G$ containing $g$ vertices. Then:
\begin{enumerate}
\item $L$ is contained in $N=\# (g+1)<\nns S^2\times S^1\nns >$;
\item $N\setminus L$ can be decomposed into the union of $g$ copies of $D(0,0,0,0,0,0)$-blocks (in one-to one correspondence with the vertices of $G$) glued to each other along their boundary spheres and hence, in particular, its volume is $2gVol_{Oct}$, where $Vol_{Oct}$ is the hyperbolic volume of a regular ideal octahedron;
\item Each pair $(M,L')$ where $M$ is a closed $3$-manifold and $L'$ is a link in $M$, can be obtained by an integer Dehn filling on some components of a suitable fundamental shadow link $L$;
\item In the preceding point, if $M\setminus L'$ is hyperbolic with volume $V$, then one can find  $L$ so that $hV\leq sc(L)\leq kV^2$ where $h$ and $k$ are two constants not depending on the data.  
\end{enumerate}
\end{teo} 

Let from now on $L=L_1\cup \ldots\cup L_r$ be a fundamental shadow link associated to a graph $G$ with $g$ vertices and let us color the components of $L$ with non-negative real numbers $a_i$; let moreover $m_i$ be the meridian of $L_i$ oriented arbitrarily. Let $w_1,\ldots w_g$ be the vertices of $G$ and for each $w_i$ let $a_j(i), j=0\ldots 5$ be the colors associated to the $6$ strands of $L$ passing through $w_i$ disposed so that $a_j(i)$ and $a_{j+3}(i)\  (mod\ 6)$ correspond to opposite strands. 
Let us denote $(N\setminus L)_{a}$ the (possibly non-complete) hyperbolic $3$-manifold obtained by equipping $N\setminus L$ with the hyperbolic structure such that, for each $i=1,\ldots r$, the holonomy of $m_i$ is conjugated to an upper triangular matrix whose eigenvalues are $e^{\pm\pi ia_i}$.   
The following is a generalization of points $2)$ and $3)$ of Theorem \ref{teo:ioedylan}:
\begin{prop}\label{prop:kerchoff}
There exists a neighborhood $U$ of $\vec 0\in \mr^r$ such that for each $\vec a\in U$ the metric completion of $(N\setminus L)_{a}$ can be decomposed into $g$ $D$-blocks $D_1,\ldots D_g$ glued to each other along their boundary spheres and where each $D_s$ is $$D_s=D(2\pi|a_0(s)|,2\pi|a_1(s)|,2\pi|a_2(s)|,2\pi|a_3(s)|,2\pi|a_4(s)|, 2\pi|a_5(s)|)$$
\end{prop} 
\begin{prf}{1}{
If one equips $D_s$ with the structure of the thesis, when gluing two different $D$-blocks along their boundary spheres according to the combinatorial structure of $G$ and $L$, the blocks geometrically match up correctly because their boundaries are geodesic spheres having the same cone angle singularities around their $3$ intersections with $L$. This way one builds an hyperbolic structure on $N$ having cone angle singularities around the $i^{th}$ component of $L$ whose total angle is $2\pi a_i$. It is well known (see for instance \cite{Th}) that one obtains a structure having the same property by completing the hyperbolic structure of $(N\setminus L)_{a}$ .
By Mostow Rigidity the two structures are isometric.
}
\end{prf}

\section{The asymptotics of $6j$-symbols}
In this section we recall the definition of $6j$-symbols of the standard representation theory of $U_q(sl_2(\mc))$. After reviewing the main known results on the asymptotical behavior of these objects, we compute the asymptotical growth rate of the family of ``hyperbolic" $6j$-symbols and identify it as the volume of suitable hyperbolic truncated tetrahedra.

\subsection{Quantum objects}
 Let $q$ be a complex variable and, for each $n\in \mathbb{N}$ let: $$
 \{n\}=(-\sqrt{-1})(q^\frac{n}{2}-q^{-\frac{n}{2}})$$
  $$ \{n\}!= \prod_{1\leq i\leq n} \{i\}$$
$$ [n]=\frac{\{n\}}{\{1\}}$$
$$[n]!= \prod_{1\leq i\leq n} [i],\ [0]!=[1]!=1$$ 
\begin{defi}\label{def:admissibility}
We say that a triple $(b_0,b_1,b_2)$ of elements of $\frac{\mathbb{N}}{2}$ is \emph{admissible} if the following conditions are satisfied:
\begin{enumerate}
\item $|b_0+b_1|\geq b_2,\ |b_0+b_2|\geq b_1,\  |b_1+b_2|\geq b_0$;
\item $b_0+b_1+b_2 \in\mathbb{N}$.
\end{enumerate}
A $6$-tuple $(b_0,b_1,b_2,b_3,b_4,b_5)\in \frac{\mn}{2}^6$ is \emph{admissible} if each of the $3$-tuples $(b_0,b_1,b_2)$, $(b_0,b_4,b_5)$, $(b_3,b_1,b_5)$ and $(b_3,b_4,b_2)$ is admissible. 
\end{defi}
For each admissible triple of elements of $\frac{\mathbb{N}}{2}$ let $$\Delta(b_0,b_1,b_2)= \sqrt{\frac{[b_0+b_1-b_2]![b_0+b_2-b_1]![b_1+b_2-b_0]!}{[b_0+b_1+b_2+1]!}}$$
For each admissible $6$-tuple $(b_0,b_1,b_2,b_3,b_4,b_5)\in \frac{\mn}{2}^6$, let $T_0=b_0+b_1+b_2,\ T_1=b_0+b_4+b_5,\ T_2=b_3+b_1+b_5,\ T_3=b_3+b_4+b_2,\ Q_1=b_0+b_3+b1+b_4,\ Q_2=b_0+b_3+b_2+b_5$ and $Q_3=b_1+b_4+b_2+b_5$; then, we define its \emph{$6j$-symbol} as follows:
\begin{gather*}
 \setlength\arraycolsep{1pt}
\left( \begin{array}{ccc}
b_0 & b_1 &b_2\\
b_3  & b_4 & b_5\\
\end{array}\right)\ns=\Delta(b_0,b_1,b_2)\Delta(b_0,b_4,b_5)\Delta(b_3,b_1,b_5) \Delta(b_3,b_4,b_2)\times\\
\times\sum_{z\in \mathbb{N}}\nns \frac{(\sqrt{-1})^{-2(b_0+b_1+b_2+b_3+b_4+b_5)} (-1)^z [z+1]!}{[z\nns -T_0]![z\nns -T_1]![z\nns-T_2]![z\nns-T_3]![ Q_1-\nns z]![Q_2 -\nns z]![Q_3-\nns z]!}
\end{gather*}
where $z$ ranges over the integers such that all the arguments of the quantum factorials are non-negative.
Remark that the above functions are holomorphic in $q^\frac{1}{2}$; infact they are ``almost" rational functions of $q^\frac{1}{2}$ (they are not because of the square root in $\Delta(i,j,k)$).  Hence one can look at their Laurent series expansion around each point of $\mc$. 

\begin{defi}\label{def:radmissibility}
Let $\vec \theta$ be an element of $[0,1]^6$. 
\begin{enumerate}
\item We say that $\vec \theta$ is {\it $\mr$-admissible} if each of the following the 3-tuples 
$$(\theta_0,\theta_1,\theta_2), (\theta_0, \theta_4,\theta_5),(\theta_3 ,\theta_2 ,\theta_5 ),(\theta_3 ,\theta_4 ,\theta_2)$$ 
satisfies inequality $\theta_i+\theta_j> \theta_k$ and those obtained by permuting indices. 
  
\item We say that $\vec \theta$ is of {\it type Reshetikhin-Turaev} if it is $\mr$-admissible and for each of the 4 $\mr$-admissible 3-tuples $(\theta_i,\theta_j,\theta_k)$ the following is true:
\begin{itemize}
\item $\theta_i < \frac{1}{2},\ \forall i$;
\item $\theta_i+\theta_j+\theta_k<1$.
\end{itemize}
\item We say that $\vec \theta$ is of {\it hyperbolic type} if it is $\mr$-admissible and, for each of the 4 $\mr$-admissible 3-tuples $(\theta_i,\theta_j,\theta_k)$ it holds: 
\begin{itemize}
\item $1<\theta_i+\theta_j+\theta_k<2$;
\item $0<\theta_i+\theta_j-\theta_k<1$, and the same holds for each permutation of the indices.
\end{itemize}
\end{enumerate}
\end{defi}
\begin{rem}
If $\vec \theta$ is of hyperbolic-type then it cannot be Reshetikhin-Turaev and vice-versa.
\end{rem}
By multiplying an admissible $3$-tuple by a factor $r>2\in \mn$ and taking its integer part, one gets an admissible $3$-tuple for the level $r$ Reshetikhin-Turaev invariants.  On contrast, admissible $3$-tuples of hyperbolic type appear while computing through shadow state-sums the colored Jones invariants of links (see Subsection \ref{sub:jones} for more details).

In what follows we will be interested in Question \ref{que:main}. For instance, if $\vec \theta$ is of Reshetikhin-Turaev type, then for each $n\in \mn$ the $6j$-symbol corresponding to $n\vec \theta$ has no poles in $q=Exp(\frac{2\pi\sqrt{-1}}{n})$ and hence the zeroth order is the first one to be computed.
Let $\vec \theta$ be of Reshetikhin-Turaev type and let $\vec b^n=(b^n_0,b^n_1,b^n_2,b^n_3,b^n_4,b^n_5),\ n\in \mn$ be a sequence of admissible $6$-tuples of half integers such that $\lim_{n\to \infty} \frac{b^n_i}{n}=\theta_i$: 
\begin{teo}[Taylor and Woodward,\cite{TW}]\label{teo:tw}
Let $T(\vec \theta)$ be the spherical tetrahedron whose edge lengths are $2\pi\theta_i,\ i=0\ldots 5$. There exists a family of continuous functions $f_n$ defined on the space of isometry classes of non-degenerate spherical tetrahedra such that the following holds:
$$\left( \begin{array}{ccc}
b^n_0 & b^n_1 & b^n_2\\
b^n_3  & b^n_4 & b^n_5\\
\end{array}\right)|_{q=Exp(\frac{2\pi\sqrt{-1}}{n})}\sim f_n(T(\vec \theta))$$ 
where $g\sim h$ means that $\lim_{n\to \infty} |g-h|=0$ point-wise.
\end{teo}
\begin{rem}
We stated Taylor and Woodward's result in a very simplified way and not in its full generality: we strongly recommend the interested reader to refer to their original paper for a detailed statement.
\end{rem}
It is worth remarking that a similar result for the asymptotical behavior of classical $6j$-symbols and euclidean tetrahedra has been proved by J. Roberts (\cite{Ro}) using completely different techniques. 

Let us contrast Theorem \ref{teo:asymptotics} with \ref{teo:tw}: in \ref{teo:tw}, $\vec\theta$ is assumed to be of Reshetikhin-Turaev type whilst in \ref{teo:asymptotics} of hyperbolic type. In \ref{teo:tw} the leading order of the $6j$-symbol associated to $\vec{\theta}$ is $0$ and the asymptotical behavior is oscillating, whilst in \ref{teo:asymptotics} the leading order is $-1$ and the asymptotical behavior is exponential. In \ref{teo:tw} $\vec{\theta}$ is associated to the edge lengths of a spherical tetrahedron whilst in \ref{teo:asymptotics} to the dihedral angles of a hyperbolic truncated tetrahedron. Finally, as we shall briefly explain in Subsection \ref{sub:jones}, the $6j$ as in \ref{teo:tw} appear while computing Reshetikhin-Turaev invariants through shadow state-sums, whilst those of \ref{teo:asymptotics} while computing the Kashaev invariants.
\subsection{The asymptotic behavior in the hyperbolic case}
Let us start by first recalling some classical results and definitions.
For each $x\in \mathbb{R}$ let us define the Lobatchevskji function $\Lambda(x)=-\int_0^x Log(|2sin(s)|)ds$; $\Lambda(x)$ is analytic out of $\{ \pi k,\ k \in \mz\}$ and $\pi$-periodic. The Lobatchevskji function is crucial to compute the volume of an ideal hyperbolic tetrahedron:
\begin{teo}[Milnor, \cite{Mi}]\label{teo:milnor}
Let $T$ be an ideal hyperbolic tetrahedron having moduli $z_0,z_1,z_2$. Then $Vol(T)=\Lambda(arg({z}_0))+\Lambda(arg({z}_1))+\Lambda(arg({z}_2))$.
\end{teo}
 
The following lemma is the first contact point between hyperbolic geometry and quantum topology: roughly speaking it states that the ``asymptotical behavior of quantum factorials is controlled by the function $\Lambda$".

\begin{lemma}\label{lem:legarouf}
Let $\alpha\in [0,1]$ and $b^n,\ n\in \mn$ be a sequence such that $\lim_{n\to \infty} \frac{b^n}{n}=\alpha$. 
The following holds:
$$\lim_{n\to \infty} -\frac{\pi}{n} Log|ev_n \{b^n\}!|=\Lambda(\pi \alpha)$$
$$\lim_{n\to \infty} -\frac{\pi}{n} Log|ev_n ( (-1)^{b^n} \{n+b^n\}!)|=\Lambda(\pi \alpha)$$
\end{lemma}
\begin{prf}{1}{
The first statement was proved by Garoufalidis and Le (\cite{GL}, Proposition 8.2). 
The latter equality is a consequence of the former. Indeed, let us first note that
$$\lim_{n\to \infty}-\frac{\pi}{n}Log(ev_n\{n+b^n\}!)=\lim_{n\to \infty}-\frac{\pi}{n}Log(ev_n\frac{\{n+b^n\}!}{\{n\}})$$
Moreover, $ev_n({n+j})=-ev_n({j}), \forall j\leq n$ and so the above limit equals:
$$\lim_{n\to \infty}-\frac{\pi}{n}Log((-1)^{b^n}ev_n(\{n-1\}!) ev_n(\{b^n\}!))$$
 
Then, we apply the first equality both to $\{n-1\}!$ and to $\{b^n\}!$ and we conclude since $\Lambda(x)$ is $\pi$-periodic and $\Lambda(0)=\Lambda(\pi)=0$.
}
\end{prf}

The following result gives a first taste of how the asymptotical behavior of quantum objects is related to geometrical ones. Let $(\theta_0, \theta_1,\theta_2)\in [0,1]^3$ be an admissible triple satisfying the condition of point 3) of Definition \ref{def:admissibility} and, without loss of generality, let us suppose $\theta_0\leq \theta_1\leq \theta_2$; let also $(b^n_0,b^n_1,b^n_2),\ n\in \mn$ be a sequence of admissible 3-tuples of half integers such that $\lim_{n\to \infty} \frac{b^n_i}{n}=\theta_i$. 
\begin{prop}\label{prop:asymp}
For sufficiently large $n\in \mn$ the function $(\Delta(b^n_0,b^n_1,b^n_2))^2$ has a pole of order $1$ in $q=Exp(\frac{2\pi \sqrt{-1}}{n})$. Moreover, it holds:
$$\lim_{n\to \infty} \frac{2\pi}{n}Log(|ev_n([n]\Delta^2(b^n_0,b^n_1,b^n_2)|)=2v(\theta_0,\theta_1,\theta_2)$$
Where $v(\theta_0,\theta_1,\theta_2)$ is given by:
$$v(\theta_0,\theta_1,\theta_2)=\Lambda(\pi (\theta_0+\theta_1+\theta_2)) -\Lambda(\pi (\theta_0+\theta_1-\theta_2))-\Lambda(\pi (\theta_0+\theta_2-\theta_1))-\Lambda(\pi (\theta_1+\theta_2-\theta_0))$$
 \end{prop}
\begin{prf}{1}{
For $n$ large enough, the 3-tuple $\vec b^n\doteq (b^n_0,b^n_1,b^n_2)$ is admissible and $\frac{\vec b^n}{n}$ satisfies condition 3) of Definition \ref{def:admissibility}.
So, for $n$ large enough, all the three factors in the numerator of $\Delta^2(b^n_0,b^n_1,b^n_2)$ are quantum factorials whose argument is strictly less than $n$ and their evaluation in $q=Exp(\frac{2\pi\sqrt{-1}}{n})$ is non zero.
In contrast, for $n$ large enough, the denominator is a quantum factorial whose argument is in the interval $(n,2n)$ and hence is a polynomial in $q$ having a root of multiplicity $1$ in $q_n=Exp(\frac{2\pi\sqrt{-1}}{n})$.

Consequently, the function $[n]\Delta^2(b^n_0,b^n_1,b^n_2)$ has no poles in $q_n$: we will now calculate the asymptotical behavior of its evaluation there when $n\to \infty$. First of all let us note the following:
$$ [n]\Delta^2(b^n_0,b^n_1,b^n_2)=\{n\}\frac{\{ b^n_0+b^n_1-b^n_2\}!\{b^n_0+b^n_2-b^n_1\}!\{b^n_1+b^n_2-b^n_0 \}!}{\{b^n_0+b^n_1+b^n_2+1\}!}.$$
Now evaluating the r.h.s. in $q_n$ and applying Lemma \ref{lem:legarouf} once per each factor in the numerator and once for the denominator divided by $\{n\}$, we get the thesis.
}\end{prf}

We are now ready to calculate the full asymptotics of hyperbolic type $6j$-symbol. For this, let us fix some notation: let $\vec \theta\in [0,1]^6$ be a 6-tuple of hyperbolic type. Let the ``squares" of $\vec \theta$ be $Q_1=\pi(\theta_0+\theta_1+\theta_3+\theta_4)$, $Q_2=\pi(\theta_0+\theta_3+\theta_2+\theta_5)$, $Q_3=\pi(\theta_1+\theta_4+\theta_2+\theta_5)$ and the ``triangles" of $\vec \theta$ be $T_0=\pi(\theta_0+\theta_1+\theta_2)$, $T_1=\pi(\theta_0+\theta_4+\theta_5)$, $T_2=\pi(\theta_3+\theta_1+\theta_5)$ and $T_3=\pi(\theta_3+\theta_4+\theta_2)$. Let also $T=Max(T_0,T_1,T_2,T_3)$. Up to permuting the indices of $\theta$ by acting through the group of symmetries of a $6j$-symbol (the symmetries of a tetrahedron), we will suppose w.l.o.g that $Q_i\leq Q_j$ if $i<j$. 
Let now $\vec b^n=(b^n_0,b^n_1,b^n_2,b^n_3,b^n_4,b^n_5)\in \frac{\mn}{2}^6,\ n\in \mn$ be a sequence of admissible $6$-tuples of half integers such that $\lim_{n\to \infty} \frac{b^n_i}{n}=\theta_i$ and let $\vec u\in \mr^6$ be the vector whose $i^{th}$ component is $4\pi(\theta_i-\frac{1}{2})$.
The following holds:
\begin{teo}\label{teo:asymptotics}
There exists a unique $z_0\in [T,min(2\pi,Q_0)]$ satisfying :
$$g(z_0)\doteq \frac{sin(2\pi-z_0)sin(Q_1-z_0)sin(Q_2-z_0)sin(Q_3-z_0)}{sin(z_0-T_0)sin(z_0-T_1)sin(z_0-T_2)sin(z_0-T_3)}=1
$$
For sufficiently large $n$, the $6j$-symbol associated to $\vec b^n$ has a pole of order $1$ in $q_n=Exp(\frac{2\pi\sqrt{-1}}{n})$, and the following holds:
\begin{equation}\label{limite}
\lim_{n\to \infty} \frac{2\pi}{n}Log (| ev_n([n] \left( \begin{array}{ccc}
b^n_0 & b^1_n & b^n_2\\
b^n_3  & b^n_4 & b^n_5\\
\end{array}\right) )|)=Vol(D(\vec u))
\end{equation}
\end{teo}

\begin{prf}{1}{
Let us first note that by hypothesis $\pi<T_i< 2\pi,\ \forall i$ and $0< Q_j-T_i< \pi,\ \forall i,j$.
It is clear that $g(x)>0,\ \forall x\in (T, min(2\pi,Q_0))$, that $\lim_{x\to T^+} g(x)=\infty$ and $\lim_{x\to min(Q_0,2\pi)^-} g(x)=0$. Hence, since $g$ is continuous on $(T,min(Q_0,2\pi))$, there exists $z_0\in (T,min(2\pi,Q_0))$ such that $g(z_0)=1$. Moreover, it holds: $$g'(x)=-g(x)(\sum_{i=0}^{i=3} Ctg(x-T_i)+\sum_{i=1}^{i=3} Ctg(Q_i-x)+Ctg(2\pi-x))$$ and, since $Ctg(x-\alpha)+Ctg(\beta-x)>0,\ \forall x\in (\alpha,\beta)$ whenever $0<\beta-\alpha<\pi$, then $g'(x)<0,\ \forall x\in (T,min(2\pi,Q_0))$. This implies that the solution $z_0$ is indeed unique.

Let us now recall that the $6j$-symbol associated to $\vec b^n$ is the product of four $\Delta(b^n_i,b^n_j,b^n_k)$ (for suitable choices of $i,j$ and $k$) and of a sum of fractions whose numerator and denominator are quantum factorials.
Moreover, for $n$ big enough, Proposition \ref{prop:asymp} shows that each of the $\Delta's$ is the square root of a rational function having a pole of order $1$ in $q_n\doteq Exp(\frac{2\pi\sqrt{-1}}{n})$ and then the product of the four $\Delta's$ has a pole of order $2$ in $q_n$. The same proposition proves that the asymptotical behavior of the product of the four  $\Delta's$ is given by the summand $V$ in the r.h.s of the formula in the statement.
To simplify the notation we set  $R^n_1=b^n_0+b^n_3+b^n_1+b^n_4$, $R^n_2=b^n_0+b^n_3+b^n_2+b^n_5$, $R^n_3=b^n_1+b^n_4+b^n_2+b^n_5$, $U^n_0=b^n_0+b^n_1+b^n_2$, $U^n_1=b^n_0+b^n_4+b^n_5$, $U^n_2=b^n_3+b^n_1+b^n_5$, $U^n_3=b^n_3+b^n_4+b^n_2$ and let $U^n=Max(U_0,U_1,U_2,U_3)$. Observe that $b^n_i,\ R^n_i,\ U^n_i$ represent the integer counterpart respectively of the entries, the squares and the triangles of $\vec\theta$, and that $R^n_i\leq R^n_j$ and $U^n_i\leq U^n_j$ if $i<j$.
The remaining part of the $6j$-symbol, let us call it $\Sigma_n$, can be expressed as follows:
\begin{gather*}
\Sigma_n= \frac{(\sqrt{-1})^{-2(\sum_i b^n_i)}}{\{1\}}\times\\
\times \sum^{z\leq R^n_1}_{z\geq U^n}\nns \frac{(-1)^z \{z+1\}!}{\{z\nns -U^n_0\}!\{z\nns -U^n_1\}!\{z\nns-U^n_2\}!\{z\nns-U^n_3\}!\{R^n_1 -\nns z\}!\{R^n_2 -\nns z\}!\{R^n_3 -\nns z\}!}
\end{gather*}
where $z$ varies in the interval $[U^n,R^n_0]\cap \mn$. 
Since by hypothesis $\pi<T$ and $Q_i-T_j<\pi,\ \forall i,j$, for $n$ big enough, it holds $n<U^n$ and $R^n_i-U^n_j<n,\ \forall i,j$. Hence, in particular, the argument of the quantum factorials in the denominator are all contained in the open interval $(0,n)\cap \mn$. Moreover, for the same reasons, for $n$ big enough, $z+1>n,\ \forall z\in [U^n,R^n_0]$. 
This implies that, when we evaluate in $q_n=Exp(\frac{2\pi\sqrt{-1}}{n})$, we get no zeros in the quantum factorials in the denominator and at least one in $\{z+1\}!$; moreover, since by hypothesis $T<2\pi$, then for $n$ big enough, $U^n<2n$ and hence at least one of the summands has a zero of multiplicity exactly $1$ at $q_n$. 

To summarize, for $n$ large enough $\Sigma_n$ has a zero of order $1$ at $q_n$, and so, the whole $6j$-symbol has a pole of order $2-1=1$ at $q_n$.
Hence, taking into account Proposition \ref{prop:asymp} we are left to compute:
$$\lim_{n\to \infty} \frac{2\pi}{n}Log(|ev_n(\frac{\Sigma_n}{[n]})|)$$
To do it, let us first concentrate on the signs of the evaluation of the summands of $\Sigma_n$ in $q_n$: we claim that they are constant. Indeed, since the argument of the denominator of each summand of $\Sigma_n$ is less than $n$, its evaluation in $q_n$ is a positive real number. On contrast, if $n<z+1<2n$, then $ev_n((-1)^z\frac{\{z+1\}!}{[n]})= (-1)^z(-1)^{z+1-n}\{z+1-n\}!=(-1)^{(1-n)}\{z+1-n\}!$; if instead $z+1\geq 2n$ then $ev_n(\frac{\{z+1\}!}{[n]})=0$. Hence, to summarize, the sign of each of the summands composing $ev_n(\frac{\Sigma_n}{[n]})$ is $(-1)^{1-n}$. 

Since the signs are constant, in order to estimate the above limit, we will find the maximal term of the sum. This term has ratio bigger than $1$ with the two adjacent terms, it is straightforward to check that it corresponds to the solution of the equation:
$$-ev_n(\frac{[z+1][R^n_1-z][R^n_2-z][R^n_3-z]}{[z-U^n_0][z-U^n_1][z-U^n_2][z-U^n_3]})=1$$
Setting $y=\pi z/n$, this can be rewritten as:
$$\frac{Sin(2\pi-y-\frac{1}{n})Sin(\pi \frac{R^n_1}{n}-y)Sin(\pi \frac{R^n_2}{n}-y)Sin(\pi\frac{R^n_3}{n}-y)}{Sin(y-\pi\frac{U^n_0}{n})Sin(y-\pi\frac{U^n_1}{n})Sin(y-\pi\frac{U^n_2}{n})Sin(y-\pi\frac{U^n_3}{n})}=1
$$ 
Arguing as for the first equation of the proof, we can prove that there exists a unique real $y_n$ solving the above equation and belonging to the interval $(\pi \frac{U^n}{n},\pi\frac{min(R^n_0,2n-1)}{n})$.
Moreover, for $n$ big enough $|y_n-nz_0|<\frac{1}{4}$ and hence the maximal term is that corresponding the integer part of $nz_0$.

By Lemma \ref{lem:legarouf}, each summand of $\Sigma_n$, is of the form $Exp(nF(x))$ with:
\begin{equation}\label{miaformula}
F(x)=2(\Lambda(2\pi-x)\nns+\nns\Lambda(Q_1-x)\nns+\nns\Lambda(Q_2-x)\nns+\nns\Lambda(Q_3-x)\nns+\nns\Lambda(x-T_0)\nns+\nns\Lambda(x-T_1)\nns+\nns\Lambda(x-T_2)\nns+\nns\Lambda(x-T_3))
\end{equation}
where $F(x)$ is defined on $[T,min(2\pi,Q_1)]$ and attains its only maximum in $z_0$. By Lemma \ref{lem:troplim} we then conclude that $Exp(nF(z_0))$ is the dominating summand of $\Sigma_n$.
Then, summing up the contribution of the $\Delta's$ and using Proposition \ref{prop:asymp} we have that the left hand side of formula \ref{limite} equals 
\begin{equation}\label{miaformula2}
F(z_0)+v(\theta_0,\theta_1,\theta_2)+v(\theta_0,\theta_4,\theta_5)+v(\theta_3,\theta_1,\theta_5)+v(\theta_3,\theta_4,\theta_2)
\end{equation}
Note that the factor $2$ in the definition of $F(x)$ comes from the fact that in the statement of the Theorem the limit considered has a factor of $\frac{2\pi}{n}$ in front of the logarithm in contrast with the factor $\frac{\pi}{n}$ used in Lemmas \ref{lem:legarouf} and \ref{lem:troplim}; on contrast, in the case of the summands coming from the $\Delta's$ , the factor $2$ is annihilated by the square roots (compare with Proposition \ref{prop:asymp}).

Let us for the moment suppose that $\theta_i\geq \frac{1}{2}$.
Converting Formula (\ref{miaformula2}) through the identity $\Lambda(\frac{x}{2})=\frac{1}{2}Im(Li_2(Exp(\sqrt{-1}x)))$, we obtain
the double of the imaginary part of Murakamy-Yano's formula (\ref{MYformula}) applied to a tetrahedron whose angles are $\alpha_i=2\pi(\theta_i-\frac{1}{2})$ with $z=Exp(-\sqrt{-1}z_0)$. Indeed $z=Exp(-\sqrt{-1}z_0)$ is a solution of Equation (\ref{eqnz}): replacing $z=Exp(-\sqrt{-1}x)$ into it and simplifying on gets
\begin{equation}
\frac{d}{dx}(U(e^{\sqrt{-1}x},T))=\pi k,\ \ \ \ k\in \mz
\end{equation} 
and so putting $x=z_0$ solves the above equation with $k=0$ because $z_0$ is a maximum for $2Im(U(e^{-\sqrt{-1}x},T))=F(x)$.
  
To conclude, let us remark that in our case the inequalities satisfied by $\vec\theta$ imply that the angles $\alpha_i$ satisfy the hypotheses of Theorem \ref{teo:frigpetr} and so by Ushijima's Theorem \ref{teo:MY} the value of (\ref{miaformula2}) is $2Vol(T(\vec\alpha))=Vol(D(\vec u))$.  

\begin{figure}
\psfrag{0}{\large $0$}
 \centerline{\includegraphics[width=4.5cm]{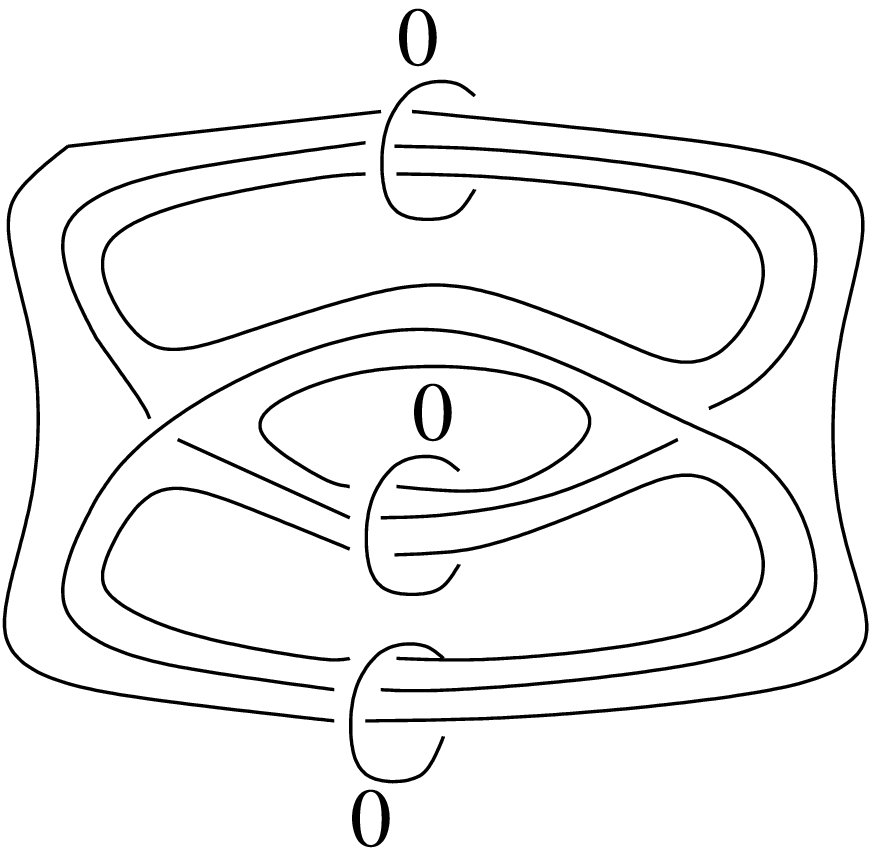}}
  \caption{}
  \label{fig:miouhl}
\end{figure}

If for some $j$ it holds $\theta_j<\frac{1}{2}$ and so $\alpha_j<0$, consider the fundamental shadow link with $6$ components depicted in Figure \ref{fig:miouhl}. It is known (\cite{NZ})  that in a neighborhood of the complete structure, the space of hyperbolic structures on the link complement is an analytic complex manifold and that, fixing a meridian $m_i,\ i=0,\ldots 5$ on each cusp one can locally parametrize it around the complete structure by considering the logarithms $u_i,\ i=0\ldots 5$ of the linear parts of the holonomies around $m_i$. Moreover, in this parametrization, the volume is a real analytic even function with respect to each $u_i$ i.e. $V(u_0,\ldots,u_i,\ldots  ,u_5)=V(u_0,\ldots,- u_i,\ldots u_5),\ \forall i$.  By Proposition \ref{prop:kerchoff}, we just proved that, when $\alpha_i$ are positive, the imaginary part of Murakamy-Yano's formula provides the volume of the hyperbolic structures such that the holonomies around the meridians of the link are conjugated to upper-triangular matrices whose eigenvalues are $Exp(\pm\sqrt{-1} \alpha_i)$. But then, being both $V$ and the imaginary part of Murakamy-Yano's formula real analytic and being identical on the set of positive $\alpha_j$, then they coincide also when $\alpha_j<1$.
}\end{prf}
\begin{rem}
The analysis in the end of the proof of the preceding theorem is similar in spirit to what is done in Section 4.2 of \cite{MY}, the main difference is that in our case the saddle point ($z_0$ in the proof) belongs to the summation interval so that saddle point method applies rigorously.  It is also important to stress that in our case $\vec\theta$ is of hyperbolic type and this, through Theorem \ref{teo:frigpetr}, allows us to associate to it a truncated tetrahedron whose angles are $\alpha_i=2\pi(\theta_i-\frac{1}{2})$; by contrast, Murakamy-Yano's analysis is inspired by the case of a proper (non truncated) tetrahedron and this would rather correspond to $\vec\theta$ of Reshetikhin-Turaev type. 
\end{rem}
\begin{lemma}\label{lem:troplim}
Let $f(x):[a,b]\to \mr^+$ be a continuous function whose maximum is attained in $x_0\in [a,b]$. For each $n\in \mn$ and each $i\in [na,nb]\cap \mn$, let $a^n_i= Exp(\frac{n}{\pi}f(i/n))$. The following holds:
$$\lim_{n\to \infty} \frac{\pi}{n}Log(\sum_{i\in[na,nb]\cap \mn} a^n_i)=f(x_0)$$ 
\end{lemma}
\begin{prf}{1}{
It is clear that $\sum a^n_i\leq n(b-a)Exp(\frac{n}{\pi}f(x_0))$ and hence that 
$$\limsup_{n\to \infty} \frac{\pi}{n}Log(\sum_{i\in[na,nb]} a^n_i)\leq f(x_0)$$
For each $\epsilon>0$ there exists a $\delta>1$ such that for each $x\in [\frac{1}{\delta}x_0,\delta x_0]$, it holds $f(x)\geq f(x_0)-\epsilon$. Now, it is evident that:
$$n(\delta-\frac{1}{\delta})Exp(\frac{n}{\pi}(f(x_0)-\epsilon))\leq \sum_{i\in[\frac{nx_0}{\delta},n\delta x_0]\cap\mn} a^n_i\leq \sum_{i\in[na,nb]\cap \mn} a^n_i$$
This easily implies that:
$$\liminf_{n\to \infty} \frac{\pi}{n}Log(\sum_{i\in[na,nb]} a^n_i)\geq f(x_0)-\epsilon$$
but since $\epsilon$ was arbitrary, the thesis follows.
}\end{prf}

\section{An application to the Volume Conjecture}\label{sec:gvc}
In this section we recall how to compute the colored Jones invariants of fundamental shadow links, this will clarify the meaning of the hyperbolic $6j$-symbols. Then, we ask a question related to a possible generalization of the Volume Conjecture for links and answer it for all the fundamental shadow links.
\subsection{Recalls on Jones invariants of fundamental shadow links}\label{sub:jones}

In \cite{Cocv} we showed how to extend the definition of colored Jones invariants to the case of links in $\# k<\nns S^2\times S^1\nns >$: in this case the resulting invariant has values in $\mq (q^{\frac{1}{2}} )$. 

To be more explicit, we remind the reader that a shadow of a link $L$ in a $3$-manifold $N$ is a $2$-dimensional spine $P$ of a $4$-manifold $W$ such that $\partial W=N$ and $P\cap \partial W=\partial P=L$.
Given a shadow $P$ of a pair $(N,L)$ as above, and a half-integer coloring $\vec b$ of the components of $L$, one can compute the Reshetikhin-Turaev invariants of the $3$-tuple $(N,L,\vec b)$ through state-sums on $P$ similar to the state-sums for Turaev-Viro invariants of spines of $3$-manifolds (see \cite{Tu}, Chapter X for a general reference). For each integer $r>2$ one can define a set of ``admissible states" (colorings of the $2$-dimensional strata of $P$) extending the coloring of $\partial P=L$ and assign a complex weight to each such state by evaluating in $Exp(2\pi i/r)$ a suitable rational function. The Reshetikhin-Turaev invariant of $(N,L,\vec b)$ of order $r$ is the sum of these weights taken over all the admissible states. 

In \cite{Cocv}, we showed that, when $N$ is $S^3$ or a connected sum of $S^2\times S^1$, for $r$ big enough, the set of admissible states stabilizes to a finite set so that the Reshetikhin-Turaev invariants can be seen as evaluations of a single rational function depending only the three-tuple $(N,L,\vec b)$; we called this rational function the colored Jones invariant of the link and denoted it $J_{\vec b}$. In the particular case when $N=S^3$ and each component of $L$ is colored by $\frac{n-1}{2}$, our colored Jones invariant equals the standard colored non-normalized Jones polynomial.
This explains why hyperbolic $6$-tuples are relevant: while expressing the $n^{th}$-colored Jones polynomial of links through shadow state-sums, if one divides the admissible states by $\frac{n-1}{2}$ one gets, near each vertex of the shadow, a $6$-tuple of hyperbolic type (see Definition \ref{def:radmissibility}).  
It is worth remarking that our arguments cannot be generalized to more general $N$ because the set of admissible states for Reshetikhin-Turaev invariants does not stabilize for $r$ big enough.

In the particular case of colored fundamental shadow links, one can give an explicit formula for the Jones invariants as follows:
\begin{prop}\label{prop:coloredjones}
Let $L=L_1\cup\ldots\cup L_r$ be a fundamental shadow link associated to a graph $G$ with $g$ vertices, and let $\vec b=b^i,\ i=1,\ldots r$ be half integer colors fixed on the components of $L$.
The colored Jones polynomial of $L$ associated to the coloring $\vec b$ is the following:
$$J_{\vec b}=\prod_{i=1}^{i=g} \left( \begin{array}{ccc}
b_0(i) & b_1(i) & b_2(i)\\
b_3(i)  & b_4(i) & b_5(i)\\
\end{array} \right) $$
where $i$ ranges over all the vertices of $G$, $b_{*}(i)$ are the six colors of the strands of $L$ passing near the $i^{th}$ vertex (so that $b_j(i)$ and $b_{j+3}(i)$ correspond to opposite strands), and a $6j$-symbol is $0$ if its entries are not admissible in the sense of Definition \ref{def:admissibility}.
\end{prop}
\begin{prf}{1}{
We limit ourselves to sketch the key-idea of the proof: we refer to \cite{Cocv} for the details. If $L$ is a fundamental shadow link, one can find a shadow $P$ of $L$ in which each $2$-dimensional stratum contains a component of $\partial P=L$. Moreover, $P$ is a simple polyhedron and contains a codimension $2$ singularity per each vertex of the graph $G$ used to describe $L$ (see Definition \ref{def:fundshlink}). This forces the state-sum to be pretty simple since there is only one admissible state which extends the coloring of the boundary; the state-sum has only one addedum which is the product of $6j$-symbols (one per each vertex of $G$).
}\end{prf}
\begin{rem}
In the above discussion we did not mention the fact that links need to be framed in order for the Jones invariants to be defined. We did not because  of the following two reasons:
\begin{itemize}
\item Changing the framing of a link changes its colored Jones invariants by multiplying it by a power of the variable, which, for the purposes of the study of the volume conjecture (and its versions) is not relevant.
\item If $L$ is a fundamental shadow link, it has a canonical framing on it such that its Jones polynomial is exactly that of Proposition \ref{prop:coloredjones}.
\end{itemize}
\end{rem}
   
\subsection{The Generalized Volume Conjecture}
In \cite{Gu}, S. Gukov proposed the following generalization of the Volume Conjecture for hyperbolic knots in $S^3$:
\begin{conj}
Let $k$ be a hyperbolic knot in $S^3$ and let $J'_{\frac{n-1}{2}}(q)$ be the $n^{th}$-colored Jones polynomial of $k$ normalized so that its value for the unknot is $1$ ($[n]J'_{\frac{n-1}{2}}=J_{\frac{n-1}{2}}$). There exists a neighborhood $U$ of $0\in \mr$ such that for each $a\in U\cap ((\mr \setminus \mq)\cup 0)$ the following holds:
$$\lim_{n\to \infty} \frac{2\pi (1+a)}{n}Log(|J'_{\frac{n-1}{2}}(e^{\frac{2\pi i (1+a)}{n}})|) =Vol_a(k)$$
where $Vol_a(k)$ is the volume of the (non-complete) hyperbolic structure on $S^3-k$ such that the holonomy around the meridian of $k$ is conjugated to a matrix whose eigenvalues are $e^{\pm (\pi i a)}$.
\end{conj}

In the same paper the conjecture was proved for the Figure Eight knot; in \cite{MY} Murakamy and Yokota provided a more general conjecture and proved it for the case of torus knots; no other proofs are known to us at present.
Remark that the above conjecture uses the deformation parameter in order to change the point where the Jones polynomials are to be evaluated. This causes, in particular, that it is hard to imagine how to extend the above conjecture for the case of links.
To do this, we propose a different approach to the problem of deforming the evaluation and ask a question which is strictly related to Gukov's conjecture.

Let us perform a small digression about the use of $ev_n$ in the statement of the conjecture. If $f$ and $g$ are meromorphic functions, then $ev_x(fg)=ev_x(f)ev_x(g)$, hence if one restates the Volume Conjectures in terms of the evaluations $ev_n$ he does no longer have to bother about the correct normalization for the Jones polynomials.
Indeed all the normalizations used in the various forms of the conjecture are of the form $[n]^k$ (for suitable fixed $k$) and since $\lim_{n\to \infty} {n}^{-1} Log(|ev_n ([n]^k)|)=0$,  they are not relevant in the limit. As a consequence, using $ev_n$, the volume conjecture makes sense even for split links and, for these, it is implied by the volume conjecture for their single components.  It then become interesting to ask the following:
\begin{question}
Let $N$ be $S^3$ or a connected sum of copies of $S^2\times S^1$ and $L\subset N$ be a link; for each $n\geq 2$, let $J_n$ be its unnormalized $n^{th}$ colored Jones polynomial. What is the asymptotical behavior of leading order in the Laurent series expansion of $J_n (L)$ around $q_n=Exp(\frac{2\pi \sqrt{-1}}{n})$? Does it have any topological meaning?  
\end{question}
It is well known that if $L$ is a knot in $S^3$ then the order is at least $1$ for each $n$ since $J_n$ is divisible by $[n]$. More in general, one could expect the above limit to be related to the presence of essential spheres in the link complement; for instance, in \cite{Cocv}, we proved that if $L$ is a fundamental shadow link in $\#_k S^2\times S^1$, the leading order is $1-k$ for every $n$. 

We stop our digression and set the following: 
\begin{defi}
Let $L\subset N=\# k\nns<\nns S^2\times S^1\nns>, k\geq 0$ be a hyperbolic link with $r$ components $L_1,\ldots L_r$, and let $m_i,\ i=1\ldots r$ be oriented meridians of $L$ and let $\vec a\in \mr^r$. We denote by $(N\setminus L)_{a}$ the hyperbolic manifold obtained by equipping $N\setminus L$ with the structure such that the holonomy around $m_i$ is conjugated to an upper triangular matrix whose eigenvalues are $e^{\pm\pi i a_j}$.
\end{defi}
\begin{question}\label{conj:genforlinks}
Let $L$ be a hyperbolic link as above.
Is it true that there exists a neighborhood $U$ of $\vec 0\in \mr^r$ such that for every $\vec a\in U\cap \mr^r$ and for each sequence $\vec b^n\in \mn^r$, $\lim_{n\to \infty}\frac{\vec b^n}{n}=\frac{\vec 1+ \vec a}{2}$, the following holds?
$$\lim_{n\to \infty} \frac{2\pi}{n}Log(|ev_{n}(J_{\vec b^n})|) =Vol((N\setminus L)_a)$$
(In the above espression, we denoted $\vec 1 \in \mr^r$ the vector whose entries are all $1$).
\end{question}

I am indebted with Stavros Garoufalidis for the following argument showing that the answer to Question \ref{conj:genforlinks} is ``no" if $L$ is a knot in $S^3$.

The symmetry principle for the normalized colored Jones polynomial states that if $k\subset S^3$ is a knot, then for all $n>m>0$ it holds:
$$J_{n\pm m}(e^{2\pi i/n})=J_{m}(e^{2\pi i/n})$$
Then, if one fixes $a=p/m\sim 1$ and lets $N=np$, $K=1/nm$, $N'=n|m-p|$, it holds:
$$J_{N}(e^{2\pi i K})=J_{N'}(e^{2\pi iK})$$
Looking at the right hand side of the above equality, one sees that $N'K=1-a\sim 0$ and as proved by Garoufalidis and Le (Theorem 4, \cite{GL}, stated below in a simplified form), the following holds:
\begin{teo}\label{teo:MMR}Ê
For every knot in $S^3$ there exists a compact neighborhood $U$ of $0\in \mc$ such that $lim_{n\to \infty} J_n(e^{2\pi i s/n})=\frac{1}{\Delta(e^{2\pi i s})}$ uniformly on $s\in U$, where $\Delta(z)$ is the Alexander polynomial of the knot normalized so that $\Delta(z)=\Delta(z^{-1})$ and $\Delta(1)=1$.
\end{teo}
Theorem \ref{teo:MMR} is a refinement of the Melvin-Morton-Rozansky Conjecture, and its proof uses in a crucial way the integrality properties of the cyclotomic function of a knot, developed by Habiro. Together with the above identity derived from the symmetry principle, it implies that $J_N(e^{2\pi i/K})$ is not exponentially growing.
These arguments apply to the normalized Jones polynomial but, when the unnormalized form is considered, the correction is given by the polynomial factor $[N]$ which cannot influence the asymptotical behavior, so the answer to Question \ref{conj:genforlinks} is ``no" for knots in $S^3$.

On contrast, for other knotted objects (such as knots in $\#_k S^2 \times S^1$, or knotted trivalent graphs) one does not have a symmetry principle, a Melvin-Morton-Rozansky Conjecture, nor a cyclotomic function with integrality properties. Thus, the above arguments do not apply to these cases and there is no evident obstruction for a positive answer to Question \ref{conj:genforlinks}.

Indeed, the following result shows that the answer is ``yes" at least for all the fundamental shadow links (all of which live in $\#_k S^2\times S^1$, $k\geq 2$).

\begin{teo}\label{teo:gcv}
The answer to Question \ref{conj:genforlinks} is ``yes" for all the fundamental shadow links. More explicitly, let $L=L_1\cup \ldots \cup L_r$ be a fundamental shadow link in $N=\# k<\nns S^2\times S^1\nns >$; there exists a neighborhood $U$ of $0\in \mr^r$ such that for every $\vec a\in U$, and for every sequence $\vec b^n\in \mn^r$ such that $\lim_{n\to \infty}\frac{\vec b^n}{n}=\frac{\vec 1+\vec a}{2}$, the following holds:
$$\lim_{n\to \infty} \frac{2\pi}{n}Log(|ev_n(J_{\vec b^n})|)=Vol((N\setminus L)_a)$$
\end{teo}
\begin{prf}{1}{
Let $G$ be a $4$-valent graph used to construct $L$ containing $g$ vertices. Let for each $n\in \mn$, let $b^n_j(i),\   j=0,\ldots ,5\ i=1,\ldots, g$ be the colors of the six strands passing through the $i^{th}$ vertex so that $b^n_j(i)$ and $b^n_{j+3}(i)$ correspond to opposite strands around the vertex (see Figure \ref{fig:uhl}) and let also $\vec u(i) \in \mr^6,\ i=1,\ldots ,g$ be the vector whose components are $$u(i)_j=\lim_{n\to \infty} 4\pi(\frac{\vec b_j^n(i)}{n}-\frac{1}{2}),\ j=0,\ldots  5$$   Remark that since by hypothesis $\vec b^n\in \mn$, the $6$-tuples $b^n_j(i), j=0,\ldots,5\ i=1,\ldots, g$ always satisfy condition 2) of Definition \ref{def:admissibility}; moreover, for $n$ large enough and $\vec a$ sufficiently near to $\vec 0\in \mr^r$, they also satisfy condition 1).  
Now, by Proposition \ref{prop:coloredjones}, it holds:
\begin{gather*}
\lim_{n\to \infty} \frac{2\pi}{n}Log(|ev_n (J_{\vec b^n})|)=\sum_{i=1}^{g} \lim_{n\to \infty}\frac{2\pi}{n}Log(ev_n( \left( \begin{array}{ccc}
b^n_0(i) & b^n_1(i) & b^n_2(i)\\
b^n_3(i)  & b^n_4(i) & b^n_5(i)\\
\end{array} \right) ))=\\
\sum_{i=1}^{g} Vol(D(\vec u(i)))=Vol((N\setminus L)_{\vec a})
\end{gather*}
where the second equality follows from Theorem \ref{teo:asymptotics} and the last equality comes from the fact that the geometric structure of $(N\setminus L)_a$ is obtained by gluing the $D(\vec u(i))$-blocks as explained in Proposition \ref{prop:kerchoff}.
}
\end{prf}

\noindent

\end{document}